\documentclass[12pt]{amsart}
\usepackage{amssymb,latexsym,amsmath,amscd}
\theoremstyle{plain}
\newtheorem{thm}{Theorem}[section]
\newtheorem{prop}[thm]{Proposition}
\newtheorem{cor}[thm]{Corollary}
\newtheorem{lemma}[thm]{Lemma}

\theoremstyle{definition}
\newtheorem{defn}[thm]{Definition}
\newtheorem{rmk}[thm]{Remark}
\newtheorem{rmks}[thm]{Remarks}
\newtheorem{notat}[thm]{Notation}

\newtheorem{ex}[thm]{Example}

\newcommand{\ra}{\rightarrow}
\newcommand{\lra}{\longrightarrow}
\newcommand{\PP}{\mathbb{P}}

\newcommand{\ZZ}{\mathbb{Z}}

\newcommand{\OO}{\mathcal{O}}

\newcommand{\K}{\mathcal{K}}

\newcommand{\I}{\mathcal{I}}

\newcommand{\MM}{\mathfrak{m}}

\newcommand{\FF}{\mathbb{F}}
\newcommand{\GG}{\mathbb{G}}

\newcommand{\OP}[1]{\OO_{\PP^{#1}}}

\newcommand{\mmax}{{\rm max}}
\newcommand{\mmin}{{\rm min}}

\newcommand{\Lra}{\Longrightarrow}
\newcommand{\hgt}{\mbox{ht}\;}
\newcommand{\rk}{\mbox{rk}\;}

\newcommand{\Coker}{\mbox{Coker}\;}
\newcommand{\Ann}{\mbox{Ann}\;}

\title{A generalized Gaeta's Theorem}

\author{Elisa Gorla}

\address{Institut f\"ur Mathematik
\\ Universit\"at Z\"urich, \hfil\break\indent Winterthurerstrasse
190, CH-8057 Z\"urich, Switzerland}

\email{elisa.gorla @ math.unizh.ch}

\thanks{The author was partially supported by the Swiss National
Science Foundation.
Part of the research in this paper was done 
while the author was a guest at the Max Planck Institut f\"ur 
Mathematik in Bonn. The author would like to thank
the Max Planck Institute for its support and hospitality.}

\begin{document}

\maketitle

{\bf Abstract:} We generalize Gaeta's
Theorem to the family of determinantal schemes.
In other words, we show that the schemes defined by minors of a
fixed size of a matrix with polynomial entries belong to the
same G-biliaison class of a complete intersection whenever
they have maximal possible codimension, given the size of the 
matrix and of the minors that define them.

\section*{Introduction}

In this paper we study the G-biliaison class of a family of schemes,
whose saturated ideals are generated by minors of matrices with 
polynomial entries. Other families of schemes defined by minors
have been studied in the same context. The results obtained in this
paper are a natural extension of some of the results proven in \cite{kl01},
\cite{ha04u2} and \cite{go05u}. In \cite{kl01} Kleppe,
Migliore, Mir\'o-Roig, Nagel, and Peterson proved that standard
determinantal schemes are glicci, i.e. that they belong to the G-liaison class
of a complete intersection. We refer
to \cite{kr00} for the definition of standard and good determinantal
schemes. Hartshorne pointed out in
\cite{ha04u2} that the double G-links produced in \cite{kl01} can indeed
be regarded as G-biliaisons. Hence, standard determinantal schemes belong 
to the G-biliaison class of a complete intersection.
In \cite{go05u} we defined symmetric
determinantal schemes as schemes whose saturated ideal is generated by the 
minors of size $t\times t$ of an $m\times m$ symmetric matrix with polynomial 
entries, and whose codimension is maximal for the given  $t$ and $m$.
In the same paper we proved that these schemes
belong to the G-biliaison class of a complete intersection.
We recently proved in \cite{go05u2} that mixed 
ladder determinantal varieties belong to the G-biliaison class of a linear variety, 
therefore they are glicci.
Ladder determinantal varieties are defined by the ideal of $t\times t$ minors of 
a ladder of inderterminates. We call them mixed ladder determinantal varieties, 
since we allow minors of different sizes in different regions of the ladder.
The results in this paper provide us with yet another family of arithmetically 
Cohen-Macaulay schemes, for which we can produce explicit G-biliaisons that 
terminate with a complete intersection. The question that one would hope to
answer is {\em whether every arithmetically Cohen-Macaulay scheme is glicci}.
Considerable progress have been made by several authors in showing that special
families of schemes are glicci (see e.g. \cite{ca00}, \cite{ca01}, \cite{kl01},
\cite{mi02}, \cite{ha02}, \cite{ca03}, \cite{ca04u}, and \cite{hu05u}).

In this paper, we study a family of schemes that correspond to ideals of minors 
of fixed size of some matrix with polynomial entries. We call them 
{\em determinantal schemes} (see Definition~\ref{det}). In Section~1 we establish
the setup, and some preliminary results about determinantal schemes.
In Remarks~\ref{ci-det} and Lemma~\ref{det&ci}, we characterize the determinantal 
schemes which are complete intersections or arithmetically Gorenstein schemes.
In Theorem~\ref{detgci} and Proposition~\ref{subm} we relate the property of being 
locally complete intersection outside a subscheme to the height of the ideal of 
minors of size one less. Section~2 contains results about heights of ideals of minors. 
It contains material that will be used to obtain the linkage results, but it can be read
independently from the rest of the article. In this section we consider an $m\times n$ 
matrix $M$, such that the ideal $I_t(M)$ has maximal height $(m-t+1)(n-t+1)$. 
In Proposition~\ref{gooddet} we show that deleting a column of $M$ we obtain a matrix $O$
whose ideal of $t\times t$ minors $I_t(O)$ has maximal height $(m-t+1)(n-t)$. In
Theorem~\ref{del}, we show that if we apply generic invertible row operations 
to $O$ and then delete a row, we obtain a matrix $N$ whose ideal of 
$(t-1)\times (t-1)$ minors has maximal height $(m-t+1)(n-t+1)$.     
Under the same assumptions, we show that if we apply generic invertible row 
operations to $M$ and then delete one entry, we obtain a ladder $L$ whose ideal of 
$t\times t$ minors has maximal height $(m-t+1)(n-t+1)-1$ (see Corollary~\ref{ladder}).
The consequence which is relevant in terms of the liaison result is that starting from 
a determinatal scheme $X$ we can produce schemes $X'$ and $Y$ such that $X'$ is determinantal 
and both $X$ and $X'$ are generalized divisors on $Y$ (see Theorem~\ref{schdiv}).
Section~3 contains the G-biliaison results. The main result of the paper is 
Theorem~\ref{bil}, where we show that any determinantal scheme can be obtained from a linear
variety by a finite sequence of ascending elementary G-biliaisons. In particular, determinantal
schemes are glicci (Corollary~\ref{glicci}). As a consequence of a result of Huneke and Ulrich, 
we obtain that determinantal schemes are in general not licci (see Corollary~\ref{licci}).

\section{Determinantal schemes}

Let $X$ be a scheme in $\PP^r=\PP^r_K$, where $K$ is an algebraically
closed field. Let $I_X$ be the saturated homogeneous ideal
associated to $X$ in the polynomial ring
$R=K[x_0,x_1,\ldots,x_r]$. For an ideal $I\subseteq R$, we denote by $H^0_*(I)$ 
the saturation of $I$ with respect to the maximal ideal 
$\MM=(x_0, x_1, \ldots ,x_r)\subseteq R$.

Let $\I_X\subseteq \OP{r}$ be the ideal sheaf of~$X$. Let $Y$ be a 
scheme that contains $X$. We denote by $\I_{X|Y}$ the ideal sheaf 
of $X$ restricted to $Y$, i.e. the quotient sheaf $\I_X/\I_Y$.
For $i\geq 0$, we let
$H^i_*(\PP^r,\I)=\oplus_{t\in\ZZ}H^i(\PP^r,\I(t))$ denote the
$i$-th cohomology module of the sheaf $\I$ on $\PP^r$. 
We simply write $H^i_*(\I)$ when there is no ambiguity on the 
ambient space $\PP^r$.

\begin{notat} Let $I\subseteq R$ be a homogeneous ideal. We let
  $\mu(I)$ denote the cardinality of a set of minimal generators
  of $I$.
\end{notat}

In this paper we deal with homogeneous ideals in the
polynomial ring $R$.

\begin{defn}
Let $M$ be a matrix with entries in $R$ of size $m\times n$, 
where $m\leq n$. 
We say that $M$ is {\em t-homogeneous} if the minors of $M$ 
of size $s\times s$ are homogeneous polynomials for all $s\leq
t$. We say that $M$ is {\em homogeneous} if its minors of any 
size are homogeneous.
\end{defn}

We always consider t-homogeneous matrices. We study a family of schemes whose 
homogeneous saturated ideal $I_t(M)$ is generated by the $t\times t$ minors of 
a t-homogeneous matrix $M$. We
regard matrices up to invertible transformations, since 
they do not change the ideal $I_t(M)$. We always assume that the
transformations that we consider preserve the $t$-homogeneity of the matrix.

\begin{defn}\label{det}
Let $X\subset\PP^r$ be a scheme. We say that $X$ is 
{\em determinantal} if: 
\begin{enumerate}
\item there exists a t-homogeneous matrix
$M$ of size $m\times n$ with entries in $R$, such that the saturated 
ideal of $X$ is generated by the minors of size $t\times t$ of $M$, 
$I_X=I_t(M)$, and
\item $X$ has codimension $(m-t+1)(n-t+1)$.
\end{enumerate}
\end{defn}

\begin{rmk}\label{codim}
The ideal $I_t(M)$ generated by the minors
of size $t\times t$ of an $m\times n$ matrix $M$ has $$\hgt
I_t(M)\leq (m-t+1)(n-t+1).$$ This is a classical result of Eagon 
and Northcott. For a proof see Theorem~2.1 in~\cite{br88b}. 
Therefore the schemes of Definition~\ref{det} have maximal
codimension for fixed $m,n,t$.
\end{rmk}

The matrix $M$ defines a morphism of free 
$R$-modules $$\varphi: R^n \lra R^m.$$
Invertible row and column operations on $M$ correspond to 
changes of basis in the domain and codomain of $\varphi$.
The scheme $X$ is the locus where $\rk\varphi\leq t-1$. 
So it only depends on the map $\varphi$ 
and not on the matrix $M$ chosen to represent it.

In some cases, we will be interested in ideals that are generated
by a subset of the minors of $M$.

\begin{notat}\label{notn}
Let $M=(F_{ij})_{1\leq i\leq m,\; 1\leq j\leq n}$ be an $m\times
n$ matrix with entries in the polynomial ring $R$. Fix a choice of row indexes $1\leq i_1\leq
i_2\leq\ldots\leq i_t\leq m$ and of column indexes $1\leq j_1\leq
j_2\leq\ldots\leq j_t\leq n$. We denote by
$M_{i_1,\ldots,i_t;j_1,\ldots,j_t}$ the determinant of the submatrix of
$M$ consisting of the rows $i_1,\ldots,i_t$ and of the columns 
$j_1,\ldots,j_t$.
\end{notat}

\begin{rmk}
Let $L$ be the subladder of $M$ consisting of all the entries except for 
$F_{mn}$. The ideal $$I_t(L)=(M_{i_1,\ldots,i_t;j_1,\ldots,j_t}\; | \; i_t\neq m 
\; \mbox{or}\; j_t\neq n)\subseteq I_t(M)$$ has height
$$\hgt I_t(L)\leq (m-t+1)(n-t+1)-1.$$ This is a special case of
Corollary~4.7 of \cite{he92a}. 
\end{rmk}

The family of determinantal schemes contains
well-studied families of schemes, such as complete intersections 
and standard determinantal schemes.

\begin{rmks}\label{ci-det}

(i) Standard determinantal schemes are a subfamily of determinantal
schemes. In fact, a determinantal scheme is standard determinantal whenever
$t=m\leq n$, that is whenever its saturated ideal is generated by the
maximal minors of $M$.

(ii) Complete intersections are a subfamily of 
determinantal schemes, since they are a subfamily of standard 
determinantal schemes. They coincide with the determinantal 
schemes that have $t=1$ or $t=m=n$ (see also Lemma~\ref{det&ci}).

(iii) The Cohen-Macaulay type of a 
determinantal scheme as of Definition~\ref{det} is
$$\prod_{i=1}^{t-1} \frac{{n-i\choose t-1}}{{m-i\choose t-1}}$$
(see \cite{br88b}). In particular, a determinantal scheme is arithmetically 
Gorenstein if and only if $m=n$. Glicciness of arithmetically Gorenstein schemes 
is established in~\cite{ca04u}. In~\cite{kl06u} it is shown that the
determinantal arithmetically Gorenstein schemes with $t+1=m=n$ are glicci.
Theorem~\ref{bil} will imply that
an arithmetically Gorenstein determinantal scheme belongs to the
G-biliaison class of a complete intersection. 
\end{rmks}
\vskip .3cm
The ideal of minors of size $t\times t$ of a generic matrix is an example of a 
determinantal scheme in $\PP^r$ for $r=mn-1$ and for each 
$t\leq m$. 

\begin{ex}\label{indet}
For any fixed $1\leq m\leq n$, and for any choice of $t$ with
$1\leq t\leq m$, let $r=mn-1$.
Let $X\subset\PP^r$
be the determinantal scheme whose saturated
ideal is generated by the minors of
size $t\times t$ of the matrix of indeterminates
$$I_X=I_t \left[\begin{array}{cccc} 
x_{1,1} & x_{1,2} & \cdots & x_{1,n} \\
x_{2,1} & x_{2,2} & \cdots & x_{2,n} \\
\vdots & \vdots & & \vdots \\
x_{m,1} & x_{m,2} & \cdots & x_{m,n}
\end{array}\right].$$
$X$ has $codim(X)=depth(I_X)=(m-t+1)(n-t+1)$ (see Theorem~2.5
of~\cite{br88b}). Then $X$ is arithmetically Cohen-Macaulay 
and determinantal. In~\cite{go05u2} we proved that $X$ belongs to the 
G-biliaison class of a complete intersection.
\end{ex}

\begin{rmk}
Complete intersections are standard determinantal, hence determinantal 
(as observed in part (ii) of Remarks~\ref{ci-det}).
Notice that the family of determinantal schemes strictly contains 
the family of standard determinantal schemes.
For example, the schemes of Example~\ref{indet} are determinantal, but not 
standard determinantal for $2\leq t\leq m-1$. This can be checked e.g. by
comparing the number of minimal generators for the saturated ideals of 
determinantal and standard determinantal schemes.
\end{rmk}

We now establish some properties of
determinantal schemes that will be needed in the sequel. 
We use the notation of Definition~\ref{det}. 
We start with a result due to~Hochster and~Eagon
(see \cite{ho71a}). We state only a special case of their
theorem, that is sufficient for our purposes. 

\begin{thm}(Hochster, Eagon)\label{he}
Determinantal schemes are arithmetically Cohen-\newline Macaulay.
\end{thm}

In the sequel, we will also need the following theorem proven by
Herzog and Trung. In Corollary 4.10 of
\cite{he92a} they establish Cohen-Macaulayness of ladder determinantal  
ideals, but we state their result only for the family of ideals that we
are interested in.

\begin{thm}(Herzog, Trung)
Let $U=(x_{ij})$ be a matrix of indeterminates of size
$m\times n$, and let $V$ be the subladder consisting 
of the all entries of $U$ except for $x_{mn}$. Then 
$$I_t(V)=(U_{i_1,\ldots,i_t;j_1,\ldots,j_t}\; | \; i_t\neq m 
\; \mbox{or}\; j_t\neq n)$$ is a Cohen-Macaulay ideal of height
$$\hgt I_t(V)=(m-t+1)(n-t+1)-1.$$
\end{thm}

We recall that if a scheme defined by the $t\times t$
minors of a matrix of indeterminates is a complete intersection, then 
it is generated by the entries of the matrix or by its 
determinant (in the case of a square matrix). We are now going to
prove the analogous result for a t-homogeneous matrix $M$ 
whose entries are arbitrary polynomials. 
We also prove a similar result for a subset of the $t\times t$
minors of $M$. 
We start by proving an easy numerical lemma.

\begin{lemma}\label{ineq}
Let $m,n,t$ be positive integers satisfying $2\leq t\leq m-1$, $m\leq n$. 
The following inequality holds:
$$(mn-t^2)(m-1)\cdot\ldots\cdot(m-t+2)(n-1)\cdot\ldots\cdot(n-t+2)>(t!)^2.$$
\end{lemma}

\begin{proof}
Since $t\leq m-1\leq n-1$,
$$(m-1)\cdot\ldots\cdot(m-t+2)(n-1)\cdot\ldots\cdot(n-t+2)\geq [(t!)/2]^2.$$
Therefore it suffices to show that
$$mn-t^2>4.$$
But $$mn-t^2\geq m^2-(m-1)^2=2m-1>4$$ since $m\geq t+1\geq 3$.
\end{proof}

The following lemma is analogous to Lemma~1.16 of~\cite{go05u}.

\begin{lemma}\label{det&ci}
Let $M$ be a t-homogeneous matrix of size $m\times n$ 
with entries in $R$ or in $R_P$ for some prime $P$. 
Let $L$ be the subladder consisting of the all entries of 
$M$ except for $F_{mn}$.
\begin{itemize}
\item[(i)] If $M$ has no invertible entries and 
$I_t(M)$ is a complete intersection of codimension $(m-t+1)(n-t+1)$, 
then $t=1$ or $t=m=n$.
\item[(ii)] If $L$ has no invertible entries and $I_t(L)$ 
is a complete intersection of codimension $(m-t+1)(n-t+1)-1$, 
then $t=1$ or $t=m=n-1$.
\end{itemize}
\end{lemma}

\begin{proof}
(i) The minors of the $t\times t$ submatrices of $M$ are a
minimal system of generators of $I_t(M)$. 
If $I_t(M)$ is a complete intersection, then 
$$\mu(I_t(M))={m\choose t}{n\choose t}=\hgt I_t(M)=(m-t+1)(n-t+1).$$ 
Computations yield 
$$[m\cdot\ldots\cdot(m-t+2)][n\cdot\ldots\cdot(n-t+2)]
=[t\cdot\ldots\cdot 2][t\cdot\ldots\cdot 2].$$
Both sides of the equality contain the same number of terms, and $t-i\leq
m-i\leq n-i$ for all $i=0,\ldots,t-2.$ So the equality holds if
and only if $t=1$ or $t=m=n$.

(ii) For a generic matrix $M=(x_{ij})$, the minors of the $t\times t$ 
submatrices that do not involve the entry $x_{mn}$ are a minimal 
system of generators of $I_t(L)$. This follows e.g. from the observation 
that they are linearly independent. By Theorem~3.5 in~\cite{br88b}, if
we substitute $F_{ij}$ for $x_{ij}$ in a minimal system of generators
of $I_t(L)$, we obtain a minimal system of generators for $I_t(L)$ in
the case $M=(F_{ij})$ and $\hgt I_t(L)=(m-t+1)(n-t+1)-1$.
In particular, the cardinality of a minimal generating system for $I_t(L)$
is in both cases 
$$\mu(I_t(L))={m\choose t}{n\choose t}-{m-1\choose t-1}{n-1\choose t-1}.$$ 
If $I_t(L)$ is a complete intersection, then 
\begin{equation}\label{hi}
\hgt I_t(L)=
{m\choose t}{n\choose t}-{m-1\choose t-1}{n-1\choose t-1}=(m-t+1)(n-t+1)-1.
\end{equation}
It follows that
$$(mn-t^2)(m-1)\cdot\ldots\cdot(m-t+1)(n-1)\cdot\ldots\cdot(n-t+1)=
(t!)^2[(m-t+1)(n-t+1)-1]$$ 
By Lemma~\ref{ineq} we have that if $t\neq 1,m$, then the left hand side of the 
equality is greater than $(t!)^2(m-t+1)(n-t+1)$. This is a contradiction, 
so $t=1$ or $t=m$.
Moreover, if $t=m$ then (\ref{hi}) simplifies to
$${n\choose m}-{n-1\choose m-1}=n-m$$
or equivalently to
$${n-1\choose m}=\frac{(n-1)\cdot\ldots\cdot(n-m)}{m!}=n-m.$$
Therefore $m=1$ or $m=n-1$. Hence either $t=1$ and
$I_t(L)$ is generated by the entries of $L$, 
or $t=m=n-1$ and $I_t(L)$ corresponds to a hypersurface (whose equation 
is the determinant of the first $m$ columns of $M$).
\end{proof}

\begin{defn}\label{loc}
Let $X\subset\PP^r$ be a scheme. We say that $X$ is {\em generically complete 
intersection} if it is locally complete intersection at all its components. 
That is, if the localization $(I_X)_P$ is generated by an
$R_P$-regular sequence for every $P$ minimal associated prime of $I_X$.

We say that $X$ is {\em locally complete 
intersection outside a subscheme of codimension d in $\PP^r$}  
if the localization $(I_X)_P$ is generated by an $R_P$-regular sequence
for every $P\supseteq I_X$ prime of $\hgt P\leq d-1$.

We say that $X$ is {\em generically Gorenstein}, abbreviated $G_0$, if it is 
locally Gorenstein at all its components. That is, if the
localization $(I_X)_P$ is a Gorenstein ideal
for every $P$ minimal associated prime of $I_X$.
\end{defn}

\begin{rmk}
The locus of points at which a scheme fails to be locally complete intersection
is closed. Therefore, a scheme of codimension $c$ in $\PP^r$ is locally complete 
intersection outside a subscheme of codimension $c+1$ in $\PP^r$ 
if and only if it is generically complete intersection.
Both of these assumption imply that the scheme is generically Gorenstein. 
\end{rmk}

We now prove two results that relate the height of the ideal 
of $(t-1)$-minors of $M$ with local properties of the scheme 
defined by the vanishing of the $t$-minors of $M$ or $L$.
The notation is as in Definition~\ref{det}.

\begin{thm}\label{detgci}
Let $X$ be a determinantal scheme with defining
matrix $M$, $I_X=I_t(M)$. Let $c=(m-t+1)(n-t+1)$ 
be the codimension of $X$. Assume that $X$ is not a complete 
intersection, i.e. $t\neq 1$ and $t,m,n$ are not all equal.
Let $d\geq c+1$ be an integer.
Then the following are equivalent:
\begin{enumerate}
\item $X$ is locally complete intersection outside of a 
subscheme of codimension $d$ in $\PP^r$.
\item $\hgt I_{t-1}(M)\geq d$. 
\end{enumerate}
\end{thm}

\begin{proof}
(1) $\Lra$ (2): let $P\supseteq I_t(M)$ be a prime ideal of height 
$c\leq\hgt P\leq d-1.$ In order to prove (2), it suffices to show 
that $P\not\supseteq I_{t-1}(M)$.
Let $M_P$ denote the localization of $M$ at $P$. The matrix 
$M_P$ can be reduced after invertible row and column 
operations to the form $$M_P=\left[\begin{array}{cc}
I_s & 0 \\ 0 & B
\end{array}\right],$$
where $I_s$ is an identity matrix of size $s\times s$, 
$0$ represents a matrix of zeroes,
and $B$ is a matrix of size $(m-s)\times (n-s)$ that has no
invertible entries.
By assumption, $I_t(M)_P\subseteq R_P$ is a
complete intersection ideal. Since $I_t(M_P)=I_{t-s}(B)$
and $B$ has no invertible entries, it follows 
by Lemma~\ref{det&ci} that either $t-s=1$, or $t-s=m-s=n-s$.
If the latter holds, then $t=m=n$ and $X$ is a hypersurface. 
Then $t-s=1$ and $I_{t-1}(M_P)=R_P$, so 
$P\not\supseteq I_{t-1}(M)$.

(2) $\Lra$ (1): let $P\supseteq I_t(M)$ be a prime of height 
$c\leq\hgt P\leq d-1.$ The thesis is proven if we show that 
$I_t(M)$ is locally generated by a regular sequence at $P$.
Since $\hgt P<\hgt I_{t-1}(M)$, then $P\not\supseteq I_{t-1}(M)$, and 
the localization $M_P$ of $M$ at $P$
can be reduced, after invertible row and column operations, to
the form $$M_P=\left[\begin{array}{cc}
I_{t-1} & 0 \\ 0 & B
\end{array}\right],$$
where $I_{t-1}$ is an identity matrix of size $(t-1)\times
(t-1)$, $0$ represents a matrix of zeroes,
and $B$ is a matrix of size $(m-t+1)\times (n-t+1).$ Since 
$PR_P\supseteq I_t(M_P)=I_1(B)$, 
we have $$\mu(I_t(M)_P)\leq (m-t+1)(n-t+1)=c=\hgt I_t(M)_P.$$
Then $I_t(M)$ is locally generated by a regular sequence at $P$.
\end{proof}

\begin{rmk}
Assume that $X$ is not a complete intersection.
For $d=c+1$, the conclusion of Theorem~\ref{detgci} can be restated 
as: $X$ is generically complete intersection if and only if 
$\hgt I_{t-1}(M)>\hgt I_t(M)$.
\end{rmk}

The implication (2) $\Lra$ (1) of Theorem~\ref{detgci} clearly holds
true without the assumption that $X$ is not a complete
intersection. The next example shows that the assumption that $X$ is
not a complete intersection is necessary for the implication (1)
$\Lra$ (2).

\begin{ex}\label{hyper}
Let $F\in R$ be a homogeneous form and consider the $t\times t$ matrix 
$$M=\left[\begin{array}{ccccc} 
F & 0 & \ldots & \ldots & 0 \\
0 & F & 0 & \ldots & 0 \\ 
\vdots  & \ddots & \ddots & \ddots & \vdots \\
\vdots  & & \ddots & \ddots & 0 \\
0 & \ldots & \ldots & 0 & F \\
\end{array}\right].$$
Let $X\subseteq\PP^r$ be the scheme with $I_X=I_t(M)=(F^t)$. Then $X$ is a 
hypersurface, hence a complete intersection, therefore locally 
complete intersection outside any subscheme. 
However the ideal $I_{t-1}(M)=(F^{t-1})$ defines a 
hypersurface in $\PP^r$, hence $\hgt I_{t-1}(M)=1$.
\end{ex}

The following proposition gives a sufficient condition for the scheme
defined by $I_t(L)$ to be generically complete intersection.

\begin{prop}\label{subm}
Let $M=(F_{ij})$ be a $t$-homogeneous matrix of size $m\times n$. 
Let $L$ be the subladder of $M$ 
consisting of all the entries except for $F_{mn}$.
Let $N$ be the submatrix obtained from $M$ by deleting the last row
and column, and let $I_{t-1}(N)$ be the ideal generated by the minors of size 
$(t-1)\times (t-1)$ of $N$.
Let $Y$ be the scheme corresponding to the ideal
$I_t(L)$. Assume that $\hgt I_t(L)=c-1=(m-t+1)(n-t+1)-1$ and
$\hgt I_{t-1}(N)=c.$ Then $Y$ is generically complete intersection.
\end{prop}

\begin{proof}
Let $P$ be a minimal associated prime of $I_Y=I_t(L)$, then $P\not\supseteq I_{t-1}(N)$.
Denote by $L_P,N_P$ the localizations of $L,N$ at $P$. Then
$N_P\subseteq L_P$ contains an 
invertible minor of size $t-1$. We can assume without loss of
generality that the minor involves the first $t-1$ rows and columns. 
After invertible row and column operations (that involve only the first
$t-1$ rows and columns) we have $$L_P=\left[\begin{array}{cc}
I_{t-1} & 0 \\ 0 & B
\end{array}\right],$$
where $B$ is the localization at $P$ of the ladder obtained by
removing the entry in the lower right corner from the submatrix of $M$
consisting of the last $m-t+1$ rows and $n-t+1$ columns. We have
$$\mu((I_Y)_P)=\mu(I_1(B))\leq (m-t+1)(n-t+1)-1=\hgt (I_Y)_P.$$
Then $I_Y$ is locally generated by a regular sequence at $P$, i.e.
$Y$ is generically complete intersection.
\end{proof}

\begin{rmk}
By Proposition~\ref{subm}, the condition that $I_{t-1}(N)=c$ implies 
that $Y$ contains a determinantal subscheme $X'$ of codimension 1, 
whose defining ideal is
$I_{X'}=I_{t-1}(N).$ Notice that whenever this is the case, $Y$ is
generically complete intersection, hence it is $G_0$. Under this
assumption we have a
concept of generalized divisor on $Y$ (see 
\cite{ha04u2} about generalized divisors). Then $X'$ is a
generalized divisor on $Y$.  
Proposition~\ref{subm} proves that the existence of such
a subscheme $X'$ of codimension 1 guarantees that $Y$ is locally a complete
intersection. Notice the analogy with standard
determinantal (\cite{kl01}) and symmetric determinantal schemes 
(\cite{go05u}).
\end{rmk}

\section{Heights of ideals of minors}

In this section we study the schemes associated to the matrix
obtained from $M$ by deleting a column, or 
a column and a generalized row. We assume that the ideal $I_t(M)$ has 
maximal height according to Remark~\ref{codim}. This section can be read 
independently from the rest of the paper.

As before, let $M$ be a $t$-homogeneous matrix of size 
$m\times n$ with entries in $R$. Assume that $I_t(M)$ 
defines a determinantal scheme $X\subset\PP^r$ of 
codimension $c=(m-t+1)(n-t+1).$ We assume that $m,n,t$ are not all
equal. In fact, if $m=n=t$ then $X$ is a hypersurface and all the
results about the heights are easily verified.

\begin{defn}
Fix a matrix $O$ of size $m\times(n-1)$.
Following \cite{kr00}, we call {\em generalized row} any row 
of the matrix obtained from $O$ by applying generic invertible 
row operations. By {\em deleting a generalized row of O}
we mean that we first apply generic invertible row
operations to $O$, and then we delete a row.
\end{defn}

We start by deleting a column of $M$ and look at the scheme 
defined by the $t\times t$ minors of the remaining columns.

\begin{prop}\label{gooddet}
Let $X\subset\PP^r$ be a determinantal scheme with 
associated matrix $M$, $I_X=I_t(M)$. Let $O$ be the 
matrix obtained from $M$ by deleting a column. Then $I_t(O)$
is the saturated ideal of a determinantal scheme $Z$ of 
codimension $(m-t+1)(n-t)$. Moreover, $Z$ is locally complete 
intersection outside a subscheme of codimension 
$(m-t+1)(n-t+1)$ in $\PP^r$.
\end{prop}

\begin{proof}
From the Lemma following Theorem~2 in~\cite{br81a}
$$\hgt I_t(M)/I_t(O)\leq m-t+1.$$
Hence $\hgt I_t(O)\geq (m-t+1)(n-t+1)-(m-t+1)=(m-t+1)(n-t)$, so 
equality holds. Then $I_t(O)$ is the saturated ideal of a
determinantal scheme $Z$ of codimension $(m-t+1)(n-t)$.
Since $\hgt I_{t-1}(O)\geq\hgt I_t(M)=(m-t+1)(n-t+1)$, by 
Theorem~\ref{detgci} $Z$ is locally complete intersection 
outside a subscheme of codimension $(m-t+1)(n-t+1)$ in $\PP^r$.
\end{proof}

\begin{notat}
We let $$\varphi:\FF\lra \GG$$ 
be the morphism of free $R$-modules associated to the matrix 
$O$, $\FF=R^{n-1}$, $\GG=R^m$.
\end{notat}

Our goal is to prove that if we delete a generalized row of $O$, 
the minors of size $t-1$ of the remaining rows define a 
determinantal scheme of the same codimension as $X$. By 
the upper-semicontinuity principle, it suffices to show 
that one can apply chosen invertible row and column 
operations to $O$, then delete a row, and obtain a matrix 
whose $t-1$ minors define a determinantal scheme.

\begin{thm}\label{del}
Let $O$ be as in Proposition~\ref{gooddet}. Deleting a 
generalized row of $O$, one obtains a matrix $N$ with
$\hgt I_{t-1}(N)=(m-t+1)(n-t+1)$.
\end{thm}

\begin{proof}
If $t=m\leq n$ then $I_m(O)$ defines a good determinantal scheme, and the 
result was proven by Kreuzer, Migliore, Nagel, and Peterson 
in~\cite{kr00}. 
Assume then that $t<m\leq n$, and consider the exact sequence
associated to the morphism $\varphi$ 
$$0\lra B\lra\FF\stackrel{\varphi}{\lra}\GG\lra \Coker\varphi\lra 0.$$
Deleting a row of $O$ corresponds to a commutative diagram with 
exact rows and columns
\begin{equation}\label{delete}
\begin{array}{rrrcr}
 & & 0 & & \\
 & & \downarrow & & \\
0 & 0 & R & & \\
\downarrow & \downarrow & \downarrow & & \\
0\lra B & \lra \FF & \stackrel{\varphi}{\lra} \GG & \lra \Coker\varphi & \lra 0 \\
\downarrow & \| & \downarrow & \downarrow & \\
0\lra B' & \lra \FF & \stackrel{\varphi'}{\lra} \GG' & \lra \Coker\varphi' 
& \lra 0 \\
 & \downarrow & \downarrow & \downarrow & \\
  & 0 & 0 & 0 & \\
\end{array}
\end{equation}
where $\varphi'$ is the morphism associated to the submatrix obtained 
from $O$ after deleting a row (possibly after applying invertible row 
operations).

We first consider the case when $m<n$. 
Since $I_m(M)$ defines a standard determinantal scheme and $O$ is 
obtained from $M$ by deleting a column, then $I_m(O)$ defines a good 
determinantal scheme (see Chapter~3 of \cite{kl01}). 
By Proposition~3.2 in~\cite{kr00}, we have that $\Coker\varphi$ is 
an ideal of positive height in $R/I_m(O)$. Then there is a minimal 
generator of $\Coker\varphi$ as an $R$-module that is non zero-divisor 
modulo $I_m(O)$. Call it $f$. Denote by $s$ the multiplication map by $f$:
\begin{equation}\label{cokernels}
0\lra R/I_m(O)\stackrel{s}{\lra}\Coker\varphi\lra\Coker s\lra 0.
\end{equation}
Since $I_m(O)+(f)\subseteq\Ann_R(\Coker s)$, $\Coker s$ is supported 
on a subscheme of codimension at least $\hgt I_m(O)+1$.
We have a commutative diagram with exact rows and columns
$$\begin{array}{rrcr}
 & 0 & 0 & \\
 & \downarrow & \downarrow & \\
 & R & \lra R/I_m(O) & \lra 0 \\
 & \downarrow & \downarrow & \\
\FF & \stackrel{\varphi}{\lra} \GG & \lra \Coker\varphi & \lra 0 \\
 & \downarrow & \downarrow & \\
 & \GG' & \stackrel{\beta}{\lra} \Coker s & \lra 0 \\
 & \downarrow & \downarrow & \\
  & 0 & 0 &
\end{array}$$
Let $\pi$ denote the morphism $\GG\lra\GG'$ in the diagram above, and define 
$\varphi'=\pi\circ\varphi$. Using the snake lemma, one can check that
$$\FF\stackrel{\varphi'}{\lra}\GG'\stackrel{\beta}{\lra}\Coker s\lra 0$$
is exact. Therefore $\Coker\varphi'=\Coker s$, and by taking kernels of 
$\varphi$ and $\varphi'$ we produce a 
diagram as (\ref{delete}).

Let $P\subseteq R$ be a prime ideal, $\hgt P\leq (m-t+1)(n-t+1)-1$. 
Since $P\not\supseteq I_{t-1}(O)$, by Proposition~16.3 in~\cite{br88b}
$\mu(\Coker(\varphi_P))\leq m-t+1$. We
claim that $P\not\supseteq I_{t-1}(N)$. If $P\not\supseteq I_m(O)$,
then the claim is proven. Therefore we can assume that $P\supseteq I_m(O)$.
Localizing at $P$ the short exact sequence (\ref{cokernels})
we have that $$\mu(\Coker(\varphi'_P))=\mu(\Coker(\varphi_P))-1\leq m-t.$$
Here $\varphi_P$ and $\varphi'_P$ denote the localization at $P$
of $\varphi$ and $\varphi'$, respectively. 
Then $P\not\supseteq I_{t-1}(N)$, again by Proposition~16.3 in~\cite{br88b}. 
Therefore the claim is proven, hence $\hgt I_{t-1}(N)=c$.

Consider now the case $t<m=n$, and consider the morphism $\psi:R^m\lra
R^{m-1}$ defined by the transposed of $O$. 
We have $\hgt I_{m-2}(O)\geq\hgt I_{m-1}(M)=4>\hgt
I_{m-1}(O)=2$. The conditions of Theorem~A2.14 in~\cite{ei95}
are satisfied, hence $\Coker\psi\subseteq R/I_{m-1}(O)$ is an ideal of
positive height. One can proceed as in the previous case,
constructing an exact sequence \begin{equation}\label{cokernels2}
0\lra R/I_{m-1}(O)\stackrel{s}{\lra}\Coker\psi\lra\Coker s\lra
0.\end{equation}
This produces a commutative diagram with exact rows and columns
$$\begin{array}{rrccr}
 & 0 & & 0 & \\
 & \downarrow & & \downarrow & \\
 & R & \lra & R/I_{m-1}(O) & \lra 0 \\
 & \downarrow & & \downarrow & \\
\FF & \stackrel{\psi}{\lra} \GG & \lra & \Coker\psi & \lra 0 \\
\| & \downarrow & & \downarrow & \\
\FF & \stackrel{\psi'}{\lra} \GG' & \lra & \Coker s & \lra 0 \\
 & \downarrow & & \downarrow & \\
  & 0 & & 0 &
\end{array}$$
Let $P\subseteq R$ be a prime ideal, $\hgt P\leq (m-t+1)^2-1$. 
Since $P\not\supseteq I_{t-1}(O)$, by Proposition~16.3 in~\cite{br88b}
$\mu(\Coker(\psi_P))\leq m-t+1$. We
claim that $P\not\supseteq I_{t-1}(N)$, where $N$ is the matrix
corresponding to $\psi'$. If $P\not\supseteq I_{m-1}(O)$,
then the claim is proven. Therefore we can assume that $P\supseteq I_{m-1}(O)$.
Localizing at $P$ the short exact sequence (\ref{cokernels2})
we have that $$\mu(\Coker(\psi'_P))=\mu(\Coker(\psi_P))-1\leq m-t.$$
Here $\psi_P$ and $\psi'_P$ denote the localization at $P$
of $\psi$ and $\psi'$, respectively. 
Then $P\not\supseteq I_{t-1}(N)$, again by Proposition~16.3 in~\cite{br88b}. 
Therefore the claim is proven, hence $\hgt I_{t-1}(N)=(m-t+1)^2$.
\end{proof}

The following is a straightforward consequence of 
Proposition~\ref{gooddet} and Theorem~\ref{del}.

\begin{cor}\label{xprime}
Let $X\subset\PP^r$ be a determinantal scheme with 
associated matrix $M$, $I_X=I_t(M)$. Delete a column of $M$, 
then a generalized row, to obtain the matrix $N$. 
Then the ideal $I_{t-1}(N)$ defines a determinantal 
scheme $X'$ of the same codimension as $X$.
\end{cor}

The next corollary is obtained by repeatedly applying 
Proposition~\ref{gooddet} and Theorem~\ref{del}.

\begin{cor}\label{ci}
Let $M$ be a $t$-homogeneous matrix of size 
$m\times n$ with entries in $R$. Assume that 
$\hgt I_t(M)=(m-t+1)(n-t+1).$ Delete $t-1$ columns and $t-1$
generalized rows. The remaining entries form a regular sequence.
\end{cor}

\begin{rmk}\label{deleting}
Under the assumptions of Corollary~\ref{ci} it is clear that 
for any submatrix $H$ consisting of $n-t+1$ columns of $M$ 
$$\hgt I_1(H)\geq I_t(M)=(m-t+1)(n-t+1).$$
What we prove in Corollary~\ref{ci} is exactly that if we apply generic 
invertible row and column operations to $M$, then pick {\em any} $n-t+1$ 
columns as $H$ and delete {\em any} $t-1$ rows of $H$, the height of the 
ideal defined by the entries does not decrease. 
So after applying generic invertible row and column operations to $M$, the 
matrix has the property that the entries of any submatrix of $M$ of size 
$(m-t+1)\times (n-t+1)$ form an $R$-regular sequence. 
\end{rmk}

\begin{thm}\label{laddcanc}
Let $M$ be as above. We assume that we have applied generic
invertible row operations to $M$, and that 
$\hgt I_t(M)=(m-t+1)(n-t+1)$. Let $L$ be the ladder
obtained from $M$ by deleting the entry in position $(m,n)$.
Let $K$ be the ladder obtained from $M$ by deleting the last row and 
column, and the entry in position $(m-1,n-1)$.
Then $$\hgt I_t(L)\geq\hgt I_{t-1}(K).$$
\end{thm}

\begin{proof}
By contradiction, suppose that $h=\hgt I_t(L)<\hgt I_{t-1}(K)$. Let $P$ 
be a minimal associated prime of $I_t(L)$ of height $h$. Then 
$P\not\supseteq I_{t-1}(K)$. Denote by $K_P, L_P$ and $M_P$ the localizations 
at $P$ of $K, L$ and $M$. Since $P\not\supseteq I_{t-1}(K)$, then 
$K_P$ contains an invertible submatrix $A$ of size $(t-1)\times (t-1).$ 
Since $K_P\subseteq L_P$, $A$ is a submatrix of $L_P$ which involves neither 
the last row, 
nor the last column. Moreover, $A$ cannot involve both row $m-1$ and
column $n-1$. To fix ideas, assume that $A$ involves the first $t-1$
rows and columns. By applying invertible row and column operations to 
$M_P$, we have $$M_P\sim\left[\begin{array}{ccc}
I_{t-1} & 0  \\
0 & B_P
\end{array}\right].$$
Notice that the row and column operations can be chosen so that they only 
affect the rows and columns of $A$. Therefore 
$B_P$ is the localization at $P$ of the submatrix $B$ obtained from $M$ by 
deleting the first $t-1$ rows and columns. The same operations yield
$$L_P\sim\left[\begin{array}{ccc}
I_{t-1} & 0  \\
0 & C_P
\end{array}\right].$$
Here $C$ is obtained from $B$ by removing the entry in 
the lower right corner, and $C_P$ denotes its localization at $P$.
By Corollary~\ref{ci} the entries of $B$, hence of $C$, 
form a regular sequence in $R$. Moreover $I_1(C)\subseteq P$, since 
$P\supseteq I_t(L)$. Therefore the entries of $C_P$ form a 
regular sequence in $R_P$, and $$\hgt I_t(L)=\hgt I_t(L_P)=\hgt I_1(C_P)=c-1.$$
But this is a contradiction, since $\hgt I_{t-1}(K)\leq c-1$.
\end{proof}

\begin{cor}\label{ladder}
Let $M$ be as above. We assume that we have applied generic
invertible row operations to $M$. Let $L$ be the ladder
obtained from $M$ by deleting the entry in the lower right corner.
If $\hgt I_t(M)=(m-t+1)(n-t+1)$, then 
$$\hgt I_t(L)=(m-t+1)(n-t+1)-1.$$  
Moreover, $I_t(L)$ is generically complete intersection.
\end{cor}

\begin{proof}
Let $L_i$ be the ladder obtained from $L$ by deleting 
the last $i$ rows and columns and the entry in position $(m-i,n-i)$, 
$1\leq i\leq t-1$.
By repeatedly applying Theorem~\ref{laddcanc}, one has
\begin{equation}\label{eqL}
\hgt I_t(L)\geq\hgt I_{t-1}(L_1)\geq\ldots\geq\hgt I_1(L_{t-1})=(m-t+1)(n-t+1)-1.
\end{equation}
The last equality follows from Corollary~\ref{ci}, where we show 
that the entries of the submatrix of $M$ consisting of the last $m-t+1$
rows and the last $n-t+1$ columns form a regular sequence (see also 
Remark~\ref{deleting}). Then $\hgt I_t(L)=(m-t+1)(n-t+1)-1.$
Let $N$ be obtained from $M$ by deleting the last row and column. 
By Theorem~\ref{del} we have $\hgt I_{t-1}(N)=(m-t+1)(n-t+1).$
Then $I_t(L)$ is generically complete intersection by 
Proposition~\ref{subm}, since 
$$\hgt I_{t-1}(N)=(m-t+1)(n-t+1)>\hgt I_t(L).$$
\end{proof}

\begin{rmk}
As a consequence of Corollary~\ref{ladder}, we obtain that
\begin{equation}\label{cod1}
\hgt I_t(M)/I_t(L)\leq 1.\end{equation} 
Since we are working under the assumption that $\hgt I_t(M)=(m-t+1)(n-t+1)$, the inequality 
(\ref{cod1}) is equivalent to $\hgt I_t(L)=(m-t+1)(n-t+1)-1$.
We believe that the inequality (\ref{cod1}) holds even 
without the assumption that $\hgt I_t(M)=(m-t+1)(n-t+1)$, 
however we were not able to prove this.
\end{rmk}

Starting from a determinatal scheme $X$ we have produced 
schemes $X'$ and $Y$ such that $X'$ is determinantal and both 
$X$ and $X'$ are generalized divisors on $Y$. We summarize these results 
in the next statement.

\begin{thm}\label{schdiv}
Let $X$ be a determinantal scheme with defining
matrix $M$, $I_X=I_t(M)$. Let $c=(m-t+1)(n-t+1)$ 
be the codimension of $X\subset\PP^r$. Let $I_t(L)$ be the ideal generated by the 
minors of size $t\times t$ of $L$, where $L$ is the subladder of $M$ 
consisting of all the entries except for $F_{mn}$ 
(after applying generic invertible row operations to $M$).
Let $N$ be the submatrix obtained from $M$ by deleting the last row and column.
Let $X'$ be the determinantal scheme with $I_{X'}=I_{t-1}(N)$.
Then $I_t(L)$ is the saturated ideal of an arithmetically 
Cohen-Macaulay, generically complete intersection scheme $Y$ 
of codimension $c-1$. $Y\supseteq X,X'$, so $X$ and $X'$ are 
generalized divisors on $Y$.
\end{thm}

\section{The Theorem of Gaeta for minors of arbitrary size}

A classical theorem of Gaeta (\cite{ga53}) proves that every arithmetically 
Cohen-Macaulay codimension 2 subscheme of $\PP^r$ can be CI-linked
in a finite number of steps to a complete intersection. The result 
was reproven and stated in the language of liaison theory 
by~Peskine and~Szpiro in \cite{pe74}. In Chapter~3 of~\cite{kl01}, Gaeta's 
Theorem is regarded as a statement about standard determinantal schemes of 
codimension 2, and extended to standard determinantal schemes of arbitrary 
codimension.
With these in mind, we wish to extend the result to the larger class of 
determinantal schemes. Determinantal schemes include the standard determinantal 
ones. More precisely, the family of standard determinantal schemes coincides 
with the determinantal schemes defined by maximal minors (see
Remark~\ref{ci-det}).

The next theorem generalizes Gaeta's Theorem, Theorem~3.6 of
\cite{kl01}, and Theorem~4.1 of \cite{ha04u2}. It is the analogous of 
Theorem~2.3 of~\cite{go05u} for a matrix that is not symmetric.
A special case of Theorem~\ref{bil} for a matrix of indeterminates 
follows also from the main result in~\cite{go05u2}.

\begin{thm}\label{bil}
Any determinantal scheme in $\PP^r$ can be obtained from a linear
variety by a finite sequence of ascending elementary G-biliaisons.
\end{thm}

\begin{proof}
Let $X\subset\PP^r$ be a determinantal scheme. We
use the notation of Definition~\ref{det}. Let $M=(F_{ij})$ be
a t-homogeneous matrix whose minors of size $t\times t$ define
$X$. Apply generic invertible row operations to $M$.

Let $c$ be the codimension of $X$, $c=(m-t+1)(n-t+1)$. 
If $t=1$ or $t=m=n$ then $X$ is a complete intersection, therefore we
can perform a finite sequence of descending elementary
CI-biliaisons to a linear variety. Therefore we assume that $t\geq 2$
and that $t<m$ if $m=n$.

Let $Y$ be the scheme with associated saturated ideal 
$$I_Y=(M_{i_1,\ldots,i_t;j_1,\ldots,j_t}\; | \; i_t\neq m 
\; \mbox{or}\; j_t\neq n).$$
By Corollary~\ref{ladder} (see also Theorem~\ref{schdiv}), 
$Y$ is arithmetically Cohen-Macaulay and generically 
complete intersection. In particular, it 
satisfies the property $G_0$. The scheme $Y$ has codimension
$c-1$, and $X$ is a generalized divisor on $Y$. Therefore a biliaison on 
$Y$ is a G-biliaison, in particular it is an even G-liaison (see \cite{kl01} 
and \cite{ha04u2} for a proof). 

Let $N$ be the matrix obtained from $M$ by deleting the last row
and column. $N$ is a t-homogeneous matrix of size 
$(m-1)\times(n-1)$.
Let $X'$ be the scheme cut out by the
$(t-1)\times(t-1)$ minors of $N$. 
By Corollay~\ref{xprime} (and Theorem~\ref{schdiv}) 
$X'$ is a generalized divisor on $Y$.
We denote by $H$ a hyperplane section divisor on $Y$.
We claim that $$X\sim X'+aH \;\;\;\mbox{for some $a>0$},$$ 
where $\sim$ denotes
linear equivalence of generalized divisors on $Y$. It follows that $X$ is
obtained by an ascending elementary biliaison from $X'$.
Repeating this argument, after $t-1$ biliaisons we
reduce to the case $t=1$, when the scheme $X$ is a complete
intersection. Then we can perform descending CI-biliaisons to a
linear variety.

Let $\I_{X|Y}$, $\I_{X'|Y}$  be the ideal sheaves on $Y$ of $X$ and
$X'$. In order to prove the claim we must show that 
\begin{equation}\label{isosh} 
\I_{X|Y}\cong \I_{X'|Y}(-a)\;\;\;\mbox{for some $a>0$}.
\end{equation}

A system of generators of $I_{X|Y}=H^0_*(\I_{X|Y})=I_t(M)/I_Y$ is
given by the images in the coordinate ring of $Y$ of the $t\times
t$ minors of $M$
$$I_{X|Y}=(M_{i_1,\ldots,i_t;j_1,\ldots,j_t} \; | \; 1\leq
i_1<i_2<\ldots<i_t\leq m, 1\leq
j_1<j_2<\ldots<j_t\leq n).$$
To keep the notation simple, we denote both an element of $R$ and
its image in $R/I_Y$ with the same symbol.
By definition, the ideal of $Y$ is generated by the minors
of size $t\times t$ of $M$, 
except for those that involve both the last row and the last column.
Therefore, a minimal system of generators of
$I_{X|Y}$ is given by 
$$I_{X|Y}=(M_{i_1,\ldots,i_{t-1},m;j_1,\ldots,j_{t-1},n} \; | \; 1\leq
i_1<\ldots<i_{t-1}\leq m-1, 1\leq j_1<\ldots<j_{t-1}\leq n-1).$$

A minimal system of generators of $I_{X'|Y}=H^0_*(\I_{X'|Y})=
I_{t-1}(N)/I_Y$ is
given by the images in the coordinate ring of $Y$ of the minors
of $N$ of size $(t-1)\times (t-1)$ 
$$I_{X'|Y}=(M_{i_1,\ldots,i_{t-1};j_1,\ldots,j_{t-1}} \; | \; 1\leq
i_1<\ldots<i_{t-1}\leq m-1, 1\leq
j_1<\ldots<j_{t-1}\leq n-1).$$
Minimality of both systems of generators can be checked with a mapping cone 
argument, using the fact that the $t\times t$ minors of $M,L,N$ are minimal 
systems of generators of $I_X, I_Y, I_{X'}$ respectively.

In order to produce an isomorphism as in (\ref{isosh}), it suffices to
observe that the ratios 
\begin{equation}\label{eqk}
\frac{M_{i_1,\ldots,i_{t-1},m;j_1,\ldots,j_{t-1},n}}
{M_{i_1,\ldots,i_{t-1};j_1,\ldots,j_{t-1}}}
\end{equation}
are all equal as elements of $H^0(\K_Y(a))$, where $\K_Y$ is the
sheaf of total quotient rings of $Y$. Then the isomorphism
(\ref{isosh}) is simply given by multiplication by that element. 
Moreover, we can compute the value of $a$ as
$$deg(M_{i_1,\ldots,i_{t-1},m;j_1,\ldots,j_{t-1},n})-
deg(M_{i_1,\ldots,i_{t-1};j_1,\ldots,j_{t-1}})=deg(F_{m,n}).$$
Equality of all the ratios in (\ref{eqk}) follows if we prove that
$$M_{i_1,\ldots,i_{t-1},m;j_1,\ldots,j_{t-1},n}\cdot
M_{k_1,\ldots,k_{t-1};l_1,\ldots,l_{t-1}}-
M_{k_1,\ldots,k_{t-1},m;l_1,\ldots,l_{t-1},n}\cdot
M_{i_1,\ldots,i_{t-1};j_1,\ldots,j_{t-1}}\in
I_Y$$
for any choice of $i,j,k,l$.
This follows from Lemma~2.4 and Lemma~2.6
in~\cite{go05u}. In those two lemmas, the result is proven in 
the case $m=n$. The proof however applies with no changes to the
situation when $m\neq n$.
This completes the proof of the 
claim and of the theorem.
\end{proof}

Theorem~\ref{bil} together with standard results in liaison
theory implies that every determinantal scheme is glicci.

\begin{cor}\label{glicci}
Every determinantal scheme $X$ can be G-bilinked in $t$ steps to a 
complete intersection, whenever $X$ is defined by the minors of
size $t\times t$ of a t-homogeneous matrix.
In particular, every determinantal scheme is glicci.
\end{cor}

Finally, we wish to emphasize that determinantal schemes are in general 
not licci, i.e. they do not belong to the CI-linkage class of a complete 
intersection. This follows from the following fundamental result established 
in~\cite{hu87} (see Corollary~5.13).

\begin{thm}(Huneke, Ulrich)
Let $I\subseteq R$ be a homogeneous ideal with minimal graded free resolution
$$0\ra\bigoplus_{j=1}^{b_c}R(-n_{c_j})\ra\cdots\ra\bigoplus_{j=1}^{b_1} 
R(-n_{1_j})\ra R\ra R/I\ra 0$$
where $c=\hgt I$. If $$\mmax\{n_{c_j}\}\leq (c-1)\mmin\{n_{1_j}\}$$ then 
$R/I$ is not licci.
\end{thm}

The theorem applies e.g. to the ideals of Example~\ref{indet}, since the shifts in the
minimal free resolution of those ideals increase linearly. Therefore we have the following.

\begin{cor}\label{licci}
Let $m,n,t\in\ZZ$ such that $2\leq t\leq m\leq n$, and 
$(m-t+1)(n-t+1)\geq 3$. Let $r=mn-1$, and 
let $X\subset\PP^r$ be the determinantal scheme whose saturated
ideal is
$$I_X=I_t \left[\begin{array}{cccc} 
x_{1,1} & x_{1,2} & \cdots & x_{1,n} \\
x_{2,1} & x_{2,2} & \cdots & x_{2,n} \\
\vdots & \vdots & & \vdots \\
x_{m,1} & x_{m,2} & \cdots & x_{m,n}
\end{array}\right].$$
Then $X$ is glicci but not licci.
\end{cor}


\begin{thebibliography}{99}

\bibitem{br81a} W. Bruns, The Eisenbud-Evans generalized principal ideal 
theorem and determinantal ideals, Proc. Amer. Math. Soc. \textbf{83} 
no. 1 (1981), 19--24

\bibitem{br88b} W. Bruns and U. Vetter, Determinantal rings,
  Lecture Notes in Mathematics 1327 (1988),
Springer-Verlag, Berlin
						  
\bibitem{ca00} M. Casanellas, R. M. Mir\'o-Roig, Gorenstein liaison of curves 
in $\PP^4$. J. Algebra \textbf{230} (2000), no. 2, 656--664

\bibitem{ca01} M. Casanellas, R. M. Mir\'o-Roig, Gorenstein liaison of divisors on 
standard determinantal schemes and on rational normal scrolls,  
J. Pure Appl. Algebra  \textbf{164}  (2001),  no. 3, 325--343

\bibitem{ca03} M. Casanellas, Gorenstein liaison of 0-dimensional schemes,  
Manuscripta Math.  \textbf{111}  (2003),  no. 2, 265--275

\bibitem{ca04u} M. Casanellas, E. Drozd, R. Hartshorne, Gorenstein liaison and 
ACM sheaves, J. Reine Angew. Math. \textbf{584} (2005), 149--171

\bibitem{co93t} A. Conca, Gr\"obner bases and determinantal
  rings, Ph.D. thesis, Universit\"at-Gesamthochschule Essen (1993)

\bibitem{ei95} D. Eisenbud, Commutative algebra with a view toward
  algebraic geometry, Graduate Texts in Mathematics \textbf{150},
  Springer-Verlag, New York (1995)

\bibitem{ga53} F. Gaeta, Ricerche intorno alle variet\`a matriciali ed ai 
loro ideali, Atti del Quarto Congresso dell'Unione Matematica Italiana, 
Taormina, 1951, vol. II, 326--328, Casa Editrice Perrella, Roma (1953)

\bibitem{go05u2} E. Gorla, Mixed ladder determinantal varieties from 
two-sided ladders, to appear in J. Pure Appl. Algebra

\bibitem{go05u} E. Gorla, The G-biliaison class of symmetric determinantal 
schemes, to appear in J. Algebra

\bibitem{go79a} S. Goto, On the Gorensteinness of determinantal
loci, J. Math. Kyoto Univ. \textbf{19}  no. 2 (1979), 371--374

\bibitem{ha02} R. Hartshorne, Some examples of Gorenstein liaison in codimension 
three, Collect. Math.  \textbf{53}  (2002),  no. 1, 21--48

\bibitem{ha04u2} R. Hartshorne, Generalized Divisors and
    Biliaison, preprint (2003) available on 
http://www.arxiv.org/abs/math.AG/0301162

\bibitem{he92a} J. Herzog and N. V. Trung, Gr\"obner bases and
  multiplicity of determinantal and Pfaffian ideals, Adv. Math.,
\textbf{96}  no. 1 (1992), 1--37

\bibitem{ho71a} M. Hochster, J. Eagon, Cohen-Macaulay rings,
  invariant theory, and the generic perfection of determinantal 
loci,  Amer. J. Math. \textbf{93}  (1971), 1020--1058

\bibitem{hu87} C. Huneke, B. Ulrich, The structure of linkage,
Ann. of Math. (2)  \textbf{126}  (1987),  no. 2, 277--334

\bibitem{hu05u} C. Huneke, B. Ulrich, Liaison of monomial ideals, preprint (2005)

\bibitem{kl01} J. O. Kleppe, J. C. Migliore, R. M. Mir\'o-Roig,
U. Nagel, and C. Peterson, Gorenstein liaison, complete intersection 
liaison invariants and unobstructedness, Mem. Amer. Math. Soc. \textbf{154} 
no. 732 (2001)

\bibitem{kl06u} J. O. Kleppe, R. M. Mir\'o-Roig, 
Ideals generated by submaximal minors, work in progress

\bibitem{kr00} M. Kreuzer, J. C. Migliore, C. Peterson, and
             U. Nagel, Determinantal schemes and Buchsbaum-Rim
             sheaves, J. Pure Appl. Algebra \textbf{150} no. 2
             (2000), 155--174
	    
\bibitem{mi02} J. Migliore, U. Nagel,  Monomial ideals and the Gorenstein 
liaison class of a complete intersection.  Compositio Math.  \textbf{133}  
(2002),  no. 1, 25--36

\bibitem{pe74} C. Peskine, L. Szpiro, Liaison des vari\'et\'es 
algé\'ebriques I, Invent. Math. \textbf{26} (1974), 271--302

\end{thebibliography}
\end{document}